\documentclass[12pt,a4paper]{amsart}
\usepackage{amsfonts,amsmath,amssymb,amsxtra,latexsym}

\newtheorem{teo}{Theorem}[section]
\newtheorem{prop}[teo]{Proposition}
\newtheorem{lema}[teo]{Lemma}

\newtheorem{remark}{Remark}

\numberwithin{equation}{section}

\newcommand{\ed}{\end{document}}

\newcommand{\real}{\mathbb{R}}

\newcommand{\N}{\mathbb{N}}

\newcommand{\dsp}{\displaystyle}
\newcommand{\sen}{\,\mbox{sen}\,}

\newcommand{\fim}{\begin{flushright}$\Box$\end{flushright}}

\newcommand{\bv}{\mathbf{v}}
\newcommand{\bu}{\mathbf{u}}
\newcommand{\ba}{\mathbf{a}}

\newcommand{\di}{{\displaystyle\int}}

\renewcommand{\thefootnote}{\fnsymbol{footnote}}

\begin{document}

\title[Steady flow in curved channels]{Steady flow for incompressible fluids
in domains with unbounded curved channels}
\author[M. M. Santos]{Marcelo M. Santos}
\address{UNICAMP (State University of Campinas)\\
IMECC-Department of Mathematics\\
Rua S\'ergio Buarque de Holanda, 651 – Cidade Universit\'aria "Zeferino Vaz" – Distr. Bar\~ao Geraldo\\
Campinas - S\~ao Paulo (SP) - Brasil\\
CEP 13083-859}
\email{msantos@ime.unicamp.br}

{\let\thefootnote\relax\footnote{Accepted for publication in Bulletin of the Brazilian Mathematical Society (2015). Special Number: proceedings of HYP2014 - Fifteenth International Conference on Hyperbolic Problems: Theory, Numerics and Applications (http://www.hyp2014.impa.br/).}}



%\date{October, 3}
\maketitle
\begin{abstract}
We give an overview on the solution of the stationary Navier-Stokes equations for non
newtonian incompressible fluids established by G. Dias and M.M. Santos ({\em Steady flow for shear thickening fluids with arbitrary fluxes}, J. Differential Equations
{\bf 252} (2012), no. 6, 3873-3898)\footnotemark[1], propose a definition for domains with unbounded
curved channels which encompasses domains with an unbounded boundary, domains with nozzles, and domains with a boundary being a punctured surface, and argue on
the existence of steady flow for incompressible fluids with arbitrary fluxes in such
domains.
\end{abstract}

\footnotetext[1]{Cf. ArXiv: 1108.3595.}
%%% ----------------------------------------------------------------------
\maketitle

\section{Shear thickening fluids}
\label{section 1}
$\indent$ 
Consider the problem of solving the Navier-Stokes equations for a stationary incompressible fluid in a domain having unbounded outlets (channels), with homogeneous Dirichlet boundary condition. If
the outlets are straight cylinders (strips, in the case of dimension two) and we ask that the fluid flow converges to given {\em Poiseuille flows} 
(parallel flows) in the ends of the outlets, the problem is known as {\em Leray problem}, and it
was solved by C. Amick \cite{amick}, under the condition that the fluxes of the given Poiseuille
flows are sufficiently small. %As far as we know, the solvability of this problem with arbitraries fluxes is an open question. 
Instead of asking that the flow converges to Poiseuille flows in the ends of the outletes, O.A. Ladyzhenskaya and V.A. Solonnikov proposed and solved in \cite{ls} the
problem of giving the flux of the fluid in each outlet. They solved this problem for arbitrary
fluxes and with the outlets not being necessarily straight cylinders. Besides, their solution has 
the property that the {\em Dirichlet integral} of the velocity field of the fluid grows at most
linearly with the direction of each outlet, and they also proved that this solution is unique under some additional smallness condition. See {\em Problem 1.1} and theorems 3.1, 2.5 and 5.5 in \cite{ls}.

Let us set some notations. The fluid velocity is denoted by $\mathbf{v}$ and the symmetric part of the velocity gradient $\mathbf{\nabla v}$, by $D(\mathbf{v})$.
In this section, we consider the fluid in an open and connect set $\Omega$ in $\real^n$, $n=2,3$, with a $C^\infty$ boundary, such that $\Omega=\bigcup_{i=0}^2\Omega_i$,
where $\Omega_0$ is a bounded subset of $\real^n$ and, for $i=1,2$, in 
different cartesian coordinate system,
$
\Omega_i=\{x\equiv (x_1,x^\prime)\in\real^n;\;(-1)x_1>0,\; x^\prime
\in\Sigma_i(x_1)\},
$
with $\Sigma_i(x_1)$ being a $C^\infty$ simply
connected open set in $\real^{n-1}$, and such that, for constants $l_1,l_2$,
$0<l_1<l_2<\infty$, %independent of $x_1$, 
%they sastify
$
 \sup_{(-1)^i x_1>0}\mbox{diam}\,\Sigma_i(x_1)\le l_2
$
and $\Omega_i$ contains the cylinder
$
\begin{array}{l}
%C_{l_1}^i=
\{x\in\real^n;\;(-1)^ix_1>0\;
\mbox{e}\;|x^\prime|<l_1\}.
\end{array}
$ 
We denote by $\Sigma$, any cross section of $\Omega$ (i.e., any bounded intersection of $\Omega$ with a
$(n\!-\!1)$-dimensional plane), $\mathbf{n}$, the orthonormal vector to $\Sigma$ pointing from $\Omega_1$ toward $\Omega_2$, and by $\Omega_t$, the truncated domain 
$\Omega_0\cup\cup_{i=1}^2\{x\in\Omega_i\,;\,(-1)^i x_1<t\}$, $t>0$, and $\Omega_{i,t-1,t}=\Omega_{i,t}\smallsetminus\overline {\Omega_{i,t-1}}$. 
Also we shall use the notations: $p\ge 2$, $p'=p/(p-1)$ and
${\mathcal D}(\Omega)$ is the space of $C_0^\infty$ divergence free vector fields
defined in $\Omega$.

Because the sets $\Omega_i$ contain straight cylinders they cannot make a turn, so we say that they are not curved channels. In Section \ref{section 2} we propose a definition for domains with unbounded channels, which includes curved channels, and %sketch a proof of the existence of steady flow for incompressible fluids with arbitrary fluxes in such domains.
argue on the existence of steady flow for incompressible fluids with arbitrary fluxes in such domains.
    
In \cite{ds}, the existence theorem \cite[Theorem 3.1]{ls} was extended for 
incompressible non newtonian fluids with the viscous stress tensor %$\mathbb S$ 
giving by $|D({\mathbf{v}})|^{p-2}D({\mathbf{v}})$, in the case $p>2$.
That is, in \cite[Theorem 2.2]{ds} we prove the following theorem:

\begin{teo}
\label{main theorem in ds} 
For any $\alpha\in\mathbb{R}$ and $p\ge 2$, the problem
\begin{equation}
\label{lady-solo}
\left\{
\begin{array}{c}
\begin{array}{rll}
\mbox{div}\{|D(\mathbf{v})|^{p-2}D(\bv )\}&=&\mathbf{v}\cdot\nabla\mathbf{v}+\nabla
{\mathcal P}\;\;\;\;\mbox{in}\;\Omega\\
\nabla\cdot\mathbf{v}&=&0\;\;\;\;\mbox{in}\;\Omega\\
\mathbf{v}&=&0\;\;\;\;\mbox{on}\;\partial\Omega\\
 \int_\Sigma\mathbf{v}\cdot\mathbf{n}&=&\alpha
\end{array}
\\
 \sup_{t>0}t^{-1}\int_{\Omega_t}|\nabla\mathbf{v}|^p<\infty\,.
\end{array}
\right. 
\end{equation}
has a weak solution $(\mathbf{v},{\mathcal P})$ in 
$W_{loc}^{1,p}(\overline{\Omega})\times L_{loc}^{p}(\overline{\Omega})$.
\end{teo}

The proof extends the technique of \cite{ls}, which first solves the Navier-Stokes equations in the truncated bounded domain $\Omega_t$ and then takes the limit with respect to the truncation parameter $t$ when $t$ goes to infinity. As in \cite{sp} and \cite{ls}, and in several subsequent papers, the velocity field ${\bv}$ is sought in the form $\bv=\bu+\ba$, where $\bu$ is the new unknown with zero flux and 
$\ba\equiv\ba(x,\delta)$, $x\in\Omega$, is a constructed divergence free vector field,
depending on the parameter $\delta>0$, with flux $\alpha$ and some other important properties.  
%In \cite{ls} it is used the vector field $\ba\equiv\ba(x,\delta)$, $x\in\Omega$, constructed in \cite{sp} (cf. \cite[p. 749]{ls}) depending on the parameter $\delta>0$. 
Among these properties satisfied by $\ba$, the property that assure that,
for some constant $c$ indepedent of $\delta$, we have that
$
\int_{\Omega_t}|\bu|^{2}|\ba|^{2}
\le c\delta^2\int_{\Omega_t}|\nabla\bu|^2
$
for all $t>0$ and all divergence free vector field $\bu\in C_0^\infty(\Omega)$, 
plays an essencial role in the analysis. Indeed, using it it is possible to obtain a priori uniform estimates in $\|\bu\|_{L^2}$ for solutions $\bv=\bu+\ba$ of problem
\eqref{lady-solo}. For the non newtonian fluids considered here, i.e. $p>2$,
it turns out that the construction of the vector field $\ba$ can be quite simplified. 
Let $\tilde\ba$ be a smooth divergence free vector field, which is bounded and has bounded derivatives in $\overline{\Omega}$, vanishes on $\partial\Omega$, and has flux one, i.e. $\int_\Sigma\tilde\ba=1$ over any cross section $\Sigma$ of $\Omega$.
Then, given $\alpha\in\real$, the vector field
$\ba=\alpha\tilde\ba$ is a vector field preserving all these properties
but having flux $\alpha$ and else satisfying the following estimates, 
for some constant $c$ depending only on $\ba$, $p$ and~$\Omega$:

\ \ \ i) $\int_{\Omega_t}|\ba|^{p^\prime}|\mathbf{\mathbf\varphi}|^{p^\prime}\leq
c\, t^{(p-2)/(p-1)}\|\nabla{\mathbf\varphi}\|_{L^p(\Omega_t)}^{p'}, \ \ \forall\ t>0,
\ \ \forall \ \mathbf{\mathbf\varphi}\in{\mathcal D}(\Omega)$;

 \ \ ii) $\int_{\Omega_{i,t-1,t}}|\nabla\ba|^p\leq c,
\ \ \forall\ t\geq 1,\ \ i=1,2$;

 \ iii) $\int_{\Omega_t}|\nabla{\ba}|^p\leq c(t+1)\;,\ \ \forall\ t\geq1$.

\noindent
For the construction of the vector field $\tilde\ba$ and the proof of properties
i)-iii), see Lemma 2.1 in \cite{ds}.

In order to obtain a regular solution in $\Omega_t$, say with $\mathbf v$ in a Sobolev space $W^{2,l}(\Omega_t)$, for some positive number $l$, we consider the
Navier-Stokes equations in $\Omega_t$ with the stress tensor modified to\break 
$\frac{1}{t}+|D(\mathbf{v})|^{p-2}D(\bv )$. More precisely, we have the following
proposition \cite[Proposition 4.1]{ds}:

\begin{prop}
\label{prop.1} Let $p\ge 2$ and $\ba$ a vector field as described above, given by \cite[Lemma 2.1]{ds}. Then, for any $t>0$, the problem 
\begin{equation}
\label{truncated}
\left\{\begin{array}{l}
\mbox{div}\{\big(\frac{1}{t}+|D(\bv)|^{p-2}\big)D(\bv)\}
={\bv}\cdot\nabla{\bv}+\nabla
{\mathcal P}\;\;\;\;\mbox{in}\;\Omega_t\\
\nabla\cdot{\bv}=0\;\;\;\;\mbox{in}\;\Omega_t\\
{\bv}=\ba\;\;\;\;\mbox{on}\;\partial\Omega_t\\
% \int_\Sigma{\bu}\cdot{\bn}=0\,.
\end{array}\right. 
\end{equation}
has a weak solution $(\bv^t,{\mathcal P^t})\in W^{1,p}(\Omega_t)\times L^p(\Omega_t)$ such that, for some $l>0$, it belongs to 
$W^{2,l}(\Omega_\tau)\times W^{1,l}(\Omega_\tau)$ for any $\tau\in (0,t)$.  
\end{prop}

The weak solution in $W^{1,p}(\Omega_t)\times L^p(\Omega_t)$ of problem \ref{truncated} is obtained by the Galerkin method and the Browder-Minty method,
taking into account the following well known estimate $
\langle|x|^{p-2}x-|y|^{p-2}y,x-y\rangle\geq
c_1|x-y|^2\left(|x|^{p-2}+|y|^{p-2}\right)\ge c_2|x-y|^p\, 
%\end{equation}
$
valid for $p\ge2$, constants $c_1,c_2$, and for all $x,y\in\real^n$, and Korn's inequality \cite{neff}\footnote{In \cite{neff}, Korn's inequality is stated for dimension three. The result in dimension two can be obtained from the one in dimension three by extending the domain $U\subset\real^2$ to 
$U\times (0,1)$ and the vector field $\bv:U\to\real^2$ to 
$(\bv,0):U\times (0,1)\to\real^3$.}. Regarding the regularity 
$(\bv^t,{\mathcal P^t})\in W^{2,l}(\Omega_\tau)\times W^{1,l}(\Omega_\tau)$, for any $\tau\in (0,t)$, it is obtained by the proof of Theorem 1.2 in \cite{veiga}.
The fact that we do not have here the homogeneous Dirichlet boundary
condition $\bv=0$ in the whole boundary $\partial\Omega_t$ does not
affect the method given in \cite{veiga} to obtain $(\bv^t,{\mathcal P^t})\in W^{2,l}(\Omega_\tau)\times W^{1,l}(\Omega_\tau)$, for any $\tau\in (0,t)$, because 
$\ba=0$ in $\partial\Omega_\tau\cap\partial\Omega$ and $\partial\Omega_\tau\cap\Omega$ is interior to $\Omega_t$, for $\tau\in (0,t)$. 

With Proposition \ref{prop.1}, the next step in the proof of the Theorem 
\ref{main theorem in ds}, following \cite{ls}, is %, fixing $\tau>0$, 
 to estimate $\|\nabla(\bv^t-\ba)\|_{L^p(\Omega_\tau)}$ uniformly with respect to $t>\tau$. 
This is achieved by the energy method with the help of Korn's inequality  \cite{neff}
and Poincare's inequality (see \cite[p.56]{galdi}). Indeed, multiplying the Navier-Stokes equation by $\bu^t=\bv^t-\ba$ and integrating by parts in $\Omega_\tau$,
if we set $z(\eta)=\int_{\eta-1}^\eta y(\tau)d\tau$, $\eta>1$, $y(\tau)=\frac{1}{\tau}\left|{\nabla\bu^t}\right|_{L^2(\Omega_\tau)}^2+
\left|{\nabla\bu^t}\right|_{L^p(\Omega_\tau)}^p$, $\tau>0$, $t\ge \tau+1$, after a lengthy computation, we arrive at $z(\eta)\leq c\eta+\Psi\left(z'(\eta)\right)$ 
for a positive constant $c$ and a nice function $\Psi$. Then by a type of reverse
Gronwall lemma \cite[Lemma 2.3]{ls} (see \cite[Lemma 3.1]{ds}) it is possible to
conclude that $y(\tau)=\frac{1}{\tau}\left|{\nabla\bu^t}\right|_{L^2(\Omega_\tau)}^2+
\left|{\nabla\bu^t}\right|_{L^p(\Omega_\tau)}^p \le c_1\tau+c_2$, for some constants $c_1,c_2$ and all $\tau>0$ and $t\ge\tau+1$. From this estimate, by a diagonalization
process, weak convergence techniques and the Browder-Minty method, we obtain a solution to the problem \eqref{lady-solo}. For the details, see \cite{ds}. 

\section{Curved channels}
\label{section 2}

As we mentioned in Section \ref{section 1}, in this section we propose a definition
for domains with unbounded channels not necessarily containing straight cylinders
and give an idea how to show  
%result extends to domains containing curved lines, going to infinity, instead of straight lines.
the existence of steady flow for incompressible fluids with arbitrary fluxes in such domains. More precisely, using some concepts from
Geometry, %we believe %that our Theorem ... holds true in the following more general form:
%we can prove the following theorem:
we argue below that the following statement is true:

\bigskip

%\begin{teo}
%\label{ext result}
{\em Let $\overline{\Omega}$ be a smooth $n$-manifold with boundary %embedded 
in $\real^n$, $n=2,3$,
with a finite number of ends $\Omega_{(i)}$, $i=1,\cdots, k$, $2\le k<~\infty$, i.e., 
$\overline{\Omega}$ is diffeomorphic to a compact smooth $n$-manifold with boundary embedded in 
$\real^n$ with $k$ points removed from its boundary. Suppose that the volumes of the cut domains
$\Omega_t$ (defined below) are of order $t$. Then, given any set of real numbers 
$\alpha_i, \ i=1,\cdots, k$, such that $\alpha_1+\cdots\alpha_k=0$, the Navier-Stokes equations in
\eqref{lady-solo}, with $p>2$, and $\Omega=\overline{\Omega}-\partial\overline{\Omega}$, has a
weak solution $\bv$ in $W^{1,p}_{\mbox{\tiny\rm loc}}(\Omega)$ having flux
$\alpha_i$ in $\Omega_{(i)}$, for each $i=1,\cdots,k$, and satisfying
the Dirichlet homogeneous boundary condition $\bv|\partial\Omega=0$.
%
%Let $\Omega$ be an open set in $\real^n$ ($n=2,3$) such that its
%boundary is a smooth $(n-~1)$-manifold with a finite number of ends 
%$\Omega_{(i)}, \ i=1,\cdots, k \ (k<~\infty)$. Then, given any set of
%real numbers $\alpha_i, \ i=1,\cdots, k$, such that $\alpha_1+\cdots\alpha_k=1$,
%the system \eqref{nspl} for power-law fluids, with $p>2$, has a weak
%solution $\bv$ in $W^{1,p}_{\mbox{\tiny\rm loc}}(\Omega)$ having flux
%$\alpha_i$ in $\Omega_{(i)}$, for each $i=1,\cdots,k$.
%\end{teo} 
}

\bigskip

Next we give more details about this statement and then give an idea for its proof.  
%The precise proof we will give in a forthcoming paper.
%\smallskip

Let $\overline{\Omega}$ be a smooth $n$-manifold with boundary 
as above. Then, by definition, there is a compact
smooth $n$-manifold with boundary $\overline{\mathcal B}$ embedded in
$\real^n$, $n=2,3$, and a diffeomorphism %\break\hfill
$H:\overline{\mathcal B}-\{p_1,\cdots,p_k\}\to \overline{\Omega}$, where 
$p_1,\cdots, p_k$ are given points in $\partial\overline{\mathcal B}$
($k<\infty$). Denote 
$\mathcal B=\overline{\mathcal B}-\partial\overline{\mathcal B}$,  
${\mathcal M}=\partial{\mathcal B}-\{p_1,\cdots,p_k\}$ and 
${\mathcal S}=\partial\Omega=\partial\overline{\Omega}$, where 
$\Omega=\overline{\Omega}-\partial\overline{\Omega}$. Then
$h:=H|{\mathcal M}$ is a diffeormorphism from ${\mathcal M}$ onto $\mathcal S$, so $\mathcal S$ is a punctured $(n\!-\!1)$-manifold, or, a
$(n\!-\!1)$-manifold with a finite number of ends $\Omega_{(i)}$, $i=1,\cdots,k$, %We do not assume that $\mathcal B$ is a simply connected domain, thus, when $n=3$, $\partial{\mathcal B}$ is a compact surface (in $\real^3$) of arbitrary genus. An end $\Omega_{(i)}$ of $\Omega$ is defined by $\Omega_{(i)}=H(V_{\varepsilon_i}(p_i))$, where
where we define an end $\Omega_{(i)}$ of $\Omega$ as follows: $\Omega_{(i)}$ is the image by $H$
of the intersection of ${\mathcal B}$ with an open ball $B_\varepsilon(p_i)$ in $\real^n$ centered at $p_i$ with radius $\varepsilon_i$, sufficiently small such that ${\mathcal B}\cap B_\varepsilon(p_i)$ is a simply connected
set. We denote this intersection by $V_{\varepsilon_i}(p_i)$. Thus, 
$\Omega_{(i)}:=H(V_{\varepsilon_i}(p_i))=H({\mathcal B}\cap B_\varepsilon(p_i))$. 
%When $n=3$, $\partial{\mathcal B}$ is a compact surface (in $\real^3$). % of arbitrary genus. 
In particular, $\Omega_{(i)}$ is an open and simply connected set in $\real^n$.
Similarly, we define an end ${\mathcal S}_{(i)}$ of $\mathcal S$ as the image by $h$
of ${\mathcal M}\cap \partial V_{\varepsilon_i}(p_i)$.  When $n=3$, ${\mathcal S}_{(i)}$ is a connected smooth surface (possibly unbounded) and when $n=2$, ${\mathcal S}_{(i)}$ is the
union of two smooth curves. Notice that, since ${\mathcal M}$ does not contain
the point $p_i$, ${\mathcal M}\cap \partial V_{\varepsilon_i}(p_i)$ is
connected when $n=3$, but is the union of two disjoints pieces of a curve
when $n=2$. %In our previous terminology, an end is an outlet to infinity. 
Now we define cross sections of $\Omega_{(i)}$ and the cut domains
$\Omega_t$ of $\Omega$, for $t\ge1$. We define a cross section $\Sigma(t)\equiv \Sigma_i(t)$
%or, more precisely, $\Sigma_i(t)$, 
of $\Omega_{(i)}$, as the image
of $\mathcal B\cap \partial V_{t^{-1}\varepsilon_i}(p_i)$ by $H$. 
Notice that $V_{t^{-1}\varepsilon_i}(p_i)\subset V_{\varepsilon_i}(p_i)$,
since $t\ge1$, and $\Sigma(t)$ is a simply connected smooth
$(n\!-\!1)$-manifold in $\Omega_{(i)}$ (without boundary). %For 
When $n=3$, the boundary of a cross section $\Sigma(t)$ is a smooth simple closed curve in 
${\mathcal S}_{(i)}=\partial\Omega_{(i)}$ which turns around ${\Omega_{(i)}}$.
In particular, it is not homotopic to a point, as it is not its preimage
by $h$ in $\mathcal M$. Indeed, this preimage is a loop (i.e. a smooth
simple closed curve) in $\mathcal M$ around $p_i$, i.e. with $p_i$ in its
interior. \ (In fact, more generally, for $n=3$ we can define a cross
section of $\Omega_{(i)}$ as any $(n\!-\!1)$-manifold in $\Omega_{(i)}$
such that its boundary is the image of a loop  in $\mathcal M$ around and
sufficiently close to $p_i$, say in 
${\mathcal M}\cap\partial V_{\varepsilon_i}(p_i)$. For $n=2$ we can also
define a cross section of $\Omega_{(i)}$ as any curve in $\Omega_{(i)}$
that is the image by $H$ of an arbitrary smooth simple curve in
$V_{\varepsilon_i}(p_i)$ connecting the two components of
${\mathcal M}\cap\partial V_{\varepsilon_i}(p_i)$.)\break\hfill
Finally, regarding the {\em cut domain} $\Omega_t$ we define it 
as being the following set:
$\Omega_t=H({\mathcal B}-\cup_{i=1}^k V_{t^{-1}\varepsilon_i}(p_i))$. 
Notice that the sets $\Omega_t$ are bounded and smooth open sets in $\real^n$ (i.e. with smooth boundaries), 
they satisfy $\Omega_{t_1}\subset\Omega_{t_2}$ if $t_1<t_2$, and $\Omega=\cup_{t\ge1}\Omega_t$. 

%\smallskip

Now that we have set terminologies, we give the idea for a proof on the existence
of steady flow in the described set $\Omega$. Analogously to Section \ref{section 1},
we search a velocity $\bv$ in the form $\bv=\bu+\ba$, where $\ba$ is a given vector
field defined in $\Omega$ such that it is divergence free, $\ba|\partial\Omega=0$, it is bounded and has bounded derivatives in $\overline{\Omega}$, and 
has flux $\alpha_i$ in each end $\Omega_{(i)}$, i.e. $\int_{\Sigma_i}\ba=\alpha_i$, for $i=1,\cdots,k$. The construction of such vector field $\ba$, as seen in Section \ref{section 1},  is an important step. % in the proof outlined in Section~. 
%Once we have that, the proof of Theorem ... can be
%adapted to prove Theorem \ref{ext result}. 
%We show how to construct such a vector field for $n=3$. 
Let us show how to construct a such vector field in the case $n=3$.
%Here, we can do this in the following way.
Let $\mathcal M$ be oriented by a normal vector field $\widetilde{\mathbf{N}}$ pointing to the exterior of $\mathcal B$. Considering the class of homotopic loops around the point $p_i$, $i=1,\cdots,k$, which we denote by $[\gamma_i]$, and assuming that any loop in $\mathcal M$ is positively oriented with respect to 
$\widetilde{\mathbf{N}}$, %we assign for each $i$ a linear functional $l_i$ 
let $l_i$ be a linear functional
(defined on the space of singular 1-chains om $\mathcal M$)
such that $l_i([\gamma_j])=\alpha_i\delta_{ij}$ (where $\delta_{ij}$ is
the Kronecker delta), $i,j=1,\cdots, k$. Then by the de Rham theorem (see e.g. \cite[\S 4.17]{warner}) there exists a closed vector field (i.e. a closed $1$-form) $\mathbf{b}_i$ on $\mathcal M$ such that $l_i$ can be identified to $\mathbf{b_i}$ through the formula 
$l_i([\gamma])=\int_\gamma\mathbf{b}_i$, for any class $[\gamma]$ of a
loop $\gamma$ in $\mathcal M$. Then if we take
$\widetilde{\mathbf{b}}:=\sum_{i=1}^{k-1}\mathbf{b}_i$ and let $\mathbf{b}$ be
the pullback of $\widetilde{\mathbf{b}}$ by $h^{-1}$, we
obtain a tangent vector field $\mathbf{b}$ on $\partial\Omega$ such
that its integral on the boundary of any cross section of the outlet
$\Omega_{(i)}$ is equal to $\alpha_i$, for $i=1,\cdots,k$. Next, we 
can extend $\mathbf{b}$ to $\Omega$, first by extending it to a tubular neighborhood $V$ of $\partial\Omega$ inside $\Omega$, by setting $\mathbf{b}(y,s)=\mathbf{b}(y)+s\mathbf{N}(y)$, for $(y,s)\in V$ (i.e.
$y\in\partial\Omega$ and $s$ in some interval $(-\epsilon_y,0)\,$), where
$\mathbf{N}$ is the unit normal vector field to $\partial\Omega$
pointing to the exterior of $\Omega$. Then we extend $\mathbf{b}$ to the
entire set $\Omega$ by multiplying it by a smooth bounded function 
$\zeta\,:\,\real^n\to\real$ such that it is equal to $1$ on $V$.  
Finally, we define $\ba$ to be the curl of the vector $\zeta\mathbf{b}$. 
Then $\ba$ is divergence free and if $\Sigma_i(t)$ is a cross section of
the outlet $\Omega_{(i)}$ with a normal vector field $\mathbf{n}_i$
pointing to infinity, by Stokes theorem and the construction of 
$\ba$, we have
$$
\begin{array}{rl}
 \di_{\Sigma_i(t)}\ba\cdot\mathbf{n}_i
=& \di_{\partial\Sigma_i(t)}\zeta\mathbf{b}
= \di_{\partial\Sigma_i(t)}\mathbf{b}
= \di_{\partial V_{t^{-1}\varepsilon_i}(p_i)}\widetilde{\mathbf{b}}\\
=&\sum_{j=1}^{k-1}\di_{\partial V_{t^{-1}\varepsilon_i}(p_i)}\mathbf{b}_j
= \sum_{j=1}^{k-1}l_j([\partial V_{t^{-1}\varepsilon_i}(p_i)])\\
=&\alpha_i
\end{array}
$$
for $i=1,\cdots, k-1$. For $i=k$ this also holds true, due to the divergence theorem, the condition $\sum_{i=1}^k\alpha_i=0$ and the fact
that $\ba$ is divergence free.
Besides, since, by hypothesis, the volumes of the cut domains $\Omega_t$
are of order $t$, i.e. $\int_{\Omega_t}\le ct$ for some constant $c$, and the vector
field $\ba$ is bounded, the estimate i) in Section \ref{section 1} holds
true. Indeed, for new constants $c$, we have
$$
\begin{array}{rl}
\di_{\Omega_t}|{\mathbf\varphi}|^{p'}|\ba|^{p'} &\le c \di_{\Omega_t}|{\mathbf\varphi}|^{p'}
\le c \ \di_{\Omega_t}|\nabla {\mathbf\varphi}|^{p'}\\
&\le c \ t^{1-p'/p} \left(\di_{\Omega_t}|\nabla {\mathbf\varphi}|^{p}\right)^{p'/p}\\
&= c \ t^{(p-2)/(p-1)} \|\nabla {\mathbf\varphi}\|_{L^p(\Omega_t)}^{p'}
%&\le c t^{(p-2)/(p-1)} \left(\di_{\Omega_t}|D(\bu)|^{p}\right)^{p'/p}
\end{array}
$$ 
for all ${\mathbf\varphi}\in{\mathcal D}(\Omega)$. 
Thus, the proof for our statement stated at the beginning of this section 
can be done by following all steps in the proof of Theorem \ref{main theorem in ds}
\cite[Theorem 2.2]{ds}. %, and assuming that the
%flux in the end $\Omega_{(i)}$ of a divergence free vector field 
%$\bv$ defined in $\Omega$ is, by definition, the integral 
%$\int_{\Sigma_i}\bv\cdot\mathbf{n}_i$.

\medskip

\noindent
{\bf Remark.} {\em In the case that the compact surface $\partial{\mathcal B}$ is of genus zero, the construction of the vector field
$\widetilde{\mathbf{b}}$ above can be simplified}. Indeed, in this
case we can assume, without loss of generality, that $\mathcal B$ is
the unit ball in $\real^3$, i.e. $\partial{\mathcal B}$ is the 
sphere $S^2$, and we can take $\widetilde{\mathbf{b}}$ as the pullback by
a stereographic projection of a linear combinations of angle forms in the plane. More precisely, let $\Pi:S^2-\{p_k\}\to\real^2$ be the
stereographic projection with projection point ({\lq\lq}north pole{\rq
\rq}) $p_k$ (we can take any point $p_1,\cdots,p_k$ as the projection point) and $\omega_i$
be the 1-form $\omega_i(x,y)
=\frac{\alpha_i/2\pi}{(x-a_i)^2+(y-b_i)^2}(-(y-b_i)dy+(x-a_i)dx)$ in
$\real^2-\{\Pi(p_i)\}$,
$i=1,\cdots,k-1$, where $(a_i,b_i)=\Pi(p_i)$.
Then $\widetilde{\mathbf{b}}=\sum_{i=1}^{k-1}\Pi^\ast\omega_i$
has the required properties.

\end{document}